\documentclass[12pt, twoside, leqno]{article}

\usepackage{amsmath,amsthm}
\usepackage{amssymb}
\usepackage{enumerate}
\usepackage{graphicx}

\newtheorem{thm}{Theorem}[section]

\newtheorem{lem}[thm]{Lemma}
\newtheorem{prop}[thm]{Proposition}

\theoremstyle{definition}

\numberwithin{equation}{section}

\usepackage{pgf,tikz}
\usetikzlibrary{arrows}

\begin{document}

\baselineskip=12.1pt

%%%%%%%%%%%%%%%%

\title{ The distance spectrum of the bipartite double cover of  strongly regular graphs}

\author{S.Morteza Mirafzal\\
Department of Mathematics \\
  Lorestan University, Khorramabad, Iran\\
E-mail: mirafzal.m@lu.ac.ir\\
E-mail: smortezamirafzal@yahoo.com
}

\date{}

\maketitle

\renewcommand{\thefootnote}{}


 % \footnote{\emph{*Corresponding author.}}
\footnote{\emph{Date}: }
\renewcommand{\thefootnote}{\arabic{footnote}}

\setcounter{footnote}{0}

%----------additions

\begin{abstract} 
A  strongly  regular  graph  with parameters $(n,d,a,c)$  is a  $d$-regular graph of order $n$,  in which every pair of adjacent vertices has exactly $a$ common neighbor(s) and every pair of nonadjacent
vertices has exactly $c$ common neighbor(s).
Let $n$ be the number of vertices of the graph $G=(V,E)$. The  distance    matrix  $D=D(G)$ of   $G$ is an $n \times n $ matrix with the rows and columns indexed by $V$ 
 such that $D_{uv} = d_{G}(u, v)=d(u,v)$, where $d_{G}(u, v)$ is the distance between the vertices $u$ and $v$ in the graph $G$. In this paper,  we are  interested in determining the distance spectrum of the bipartite double cover of the family of  strongly regular graphs. In other words,  let $G=(V,E)$ be a strongly regular graph with parameters $(n,k,a,c)$. We show that there is a close relationship between the spectrum of $G$ and the distance spectrum of $B(G)$, where $B(G)$ is the double cover of $G$. We explicitly determine the distance spectrum of the graph $B(G)$,    according to the spectrum of $G$.    In fact,  according to the parameters of the graph $G$.

\end{abstract}

%\maketitle
\section{ Introduction and Preliminaries}  
In this paper, a graph $G=(V,E)$ is
considered as an undirected simple finite graph with  the vertex-set $V=V(G)$
and  the edge-set $E=E(G)$.  The
standard terminology and notation can be found in [6,8,23].\\
Let $G=(V,E)$ be a  graph. The adjacency matrix $A $ of $G$ is the square matrix with the rows and columns
indexed by the vertex-set of V such that $A_{v,w}= 1$ when $v$  is adjacent to $w$ and $A_{v,w} = 0$
otherwise. The matrix $A$ considered as a real matrix and it is clear that $A$ is symmetric. A nonzero (column) vector $u$,
indexed by $V,$  is an eigenvector of $A$ with eigenvalue $\lambda$ when $Au$ = $\lambda u$, That is,
 $$\sum_ { \  w,  \  w\leftrightarrow v }A_{v,w}u_w = \lambda u_v,$$
  for each $v \in V,$ where $w\leftrightarrow v$ means that $w$ is adjacent to $v$. In such a case, $\lambda$ is called an eigenvalue of $A$ corresponding (belonging) to the eigenvector $u$. When $\lambda$ is an eigenvalue of the matrix  $A$, then it is a zero of the  polynomial  $P(G; x)=P(x) = |xI-A|=det(xI-A)$. The polynomial  $P(x)$ is called the $characteristic$  $polynomial$  of $G$ (or the  $characteristic$  $polynomial$  of   adjacency matrix $A$).  The $geometric$ $multiplicity$ of an eigenvalue $\lambda$ is the dimension of its eigenspace.   The $algebraic$ $multiplicity$ of an eigenvalue $\lambda$ is the multiplicity of $\lambda$ as a root of the characteristic polynomial $P(x)$. Since $A$ is a real symmetric matrix, then the geometric multiplicity and algebraic multiplicity of each of its eigenvalue $\lambda$ are the same [14]. This common value is called the $multiplicity$ of $\lambda$.
The $spectrum$  of $G$  is the (multi)set of all eigenvalues of $A$  and is  denoted by $ Spec(G)=\{\lambda_1,\lambda_2,\cdots,\lambda_n\}$ and usually indexed such that  $\lambda_1\geq\lambda_2\dots\geq \lambda_n$. If the eigenvalues of $G$ are ordered by
$ \lambda_1 > \lambda_2 > \dots > \lambda_r  $, and their multiplicities are $ m_1, m_2,\dots,m_r $,  respectively,
 then we write, \\\\
 \centerline{$ Spec( G) $   =  ${ \lambda_1,\lambda_2,\dots, \lambda_r } \choose{ m_1, m_2,\dots,m_r } $ 
 or  $Spec( G)  $ = $ \{  \lambda_1^{m_1},  \lambda_2^{m_2},\dots, \lambda_r^{m_r}   \}.$ } \\\\
A graph
is called $integral$   if all of  its eigenvalues are integers. The study of integral graphs  was initiated by Harary
and Schwenk in 1974 [9].  A survey of papers up to 2002 has been  appeared in [2], but
more than a hundred new studies on integral graphs have been published in the last
 23  years (see [15,18] and references in them). \newline 
Let $n$ be the number of vertices of the graph $G$. The $distance$  $matrix$ $D=D(G)$ is an $n \times n $ matrix with the rows and columns  indexed by $V$, such that $D_{uv} = d_{G}(u, v)=d(u,v)$, where $d_{G}(u, v)$ is the distance between the vertices $u$ and $v$ in the graph $G$.
A graph $G$ is called $distance$  $integral$  (briefly, $D$-$integral$) if all of the distance eigenvalues  of $G$ are integers.
The distance matrix and distance eigenvalues of graphs have been studied by researchers for many years (see [1,7,11,13,26,27]). 
Although there are many
 papers that study distance spectrum of graphs and their applications, the $D$-integral graphs are studied only in a few papers. Some of the recent papers include [5,10,12,19,20,21,22,24,28].\\
Strongly regular graphs   are simple regular graphs with the property that the
number of common neighbors of a pair of distinct vertices depends only on whether
the two vertices are adjacent or not. They have been originally introduced by R. C. Bose
 [4,8] and they are one of the central notions of modern algebraic graph theory. Small
examples include the pentagon $C_5$, the Petersen graph, triangular graphs and  the Clebsch
graph [4,8]. Formally, a  $strongly \ regular$ graph  with parameters $(n,k,a,c)$  is a  $k$-regular graph of order $n$,  in which every pair of adjacent vertices has exactly $a$ common neighbor(s) and every pair of nonadjacent
vertices has exactly $c$ common neighbor(s). It is known and easy to check that 
the Petersen graph is a strongly regular graph with parameters $(10,3,0,1)$.\\
Let $G_1=(V_1,E_1),G_2=(V_2,E_2)  $ be  graphs. Then their direct product is the graph
$  G_1 \times G_2 $
 with the  vertex-set $ \{( v_1,v_2) \  |  \  v_1 \in V_1,   v_2 \in V_2 \} $, and for
  which vertices $( v_1,v_2)$ and $ ( w_1,w_2)  $ are adjacent precisely if $ v_1$ is adjacent to $w_1$ in $G_1$ and $ v_2$ is adjacent to $w_2$ in $ G_2$. When $G_2=K_2$, the complete graph on two vertices,  then $G \times K_2$ is known as the $bipartite$ $double$ $cover$ of the graph $G$, denoted by $B(G)$.  The notion of bipartite double cover of a graph is one the important subjects in algebraic graph theory and some of the interesting families of graphs are bipartite double covers [4].\\
Let $n \geq 3$ be an integer.  A $crown$  $graph$  $Cr(n)$ is a graph obtained from the complete bipartite graph $K_{n,n}$ by removing a perfect matching. It is easy to check that the graph $Cr(n)$ is an $(n-1)$-regular bipartite graph of diameter $3$.  The $bipartite$  $Kneser$  graph $H(n,k)$, $1 \leq k \leq n-1$, is a bipartite graph with the vertex-set consisting of all $k$-subsets and $(n-k)$-subsets of the set $[n]=\{1,2,3,\dots,n\}$, in which two vertices $v$ and $w$ are adjacent if and only if $v \subset w$ or $w \subset v$. It is easy to see that the crown graph $Cr(n)$ is isomorphic with the bipartite graph 
$H(n,1)$  [16,17].  
  Moreover,  it can be shown that the crown graph $Cr(n)$ is isomorphic with the graph $K_n \times K_2$, where $K_n$ is the complete graph on $n$ vertices [17]. Also, the bipartite Kneser graph $H(n,k)$ is isomorphic with the   bipartite double cover of the Kneser graph $K(n,k)$ [17].\\
The Clebsch graph is a strongly regular graph of parameters $(16,5,0,2)$.
In fact, it is the unique strongly regular graph with these parameters [4,8].
It can be check that the bipartite double cover of the Clebsch graph is isomorphic with the hypercube $Q_5$ [4,8,18].
Figure 1. displays a version of  the Clebsch graph in the plane. 
\begin{figure} [ht]
\centerline{\includegraphics[width=11 cm]{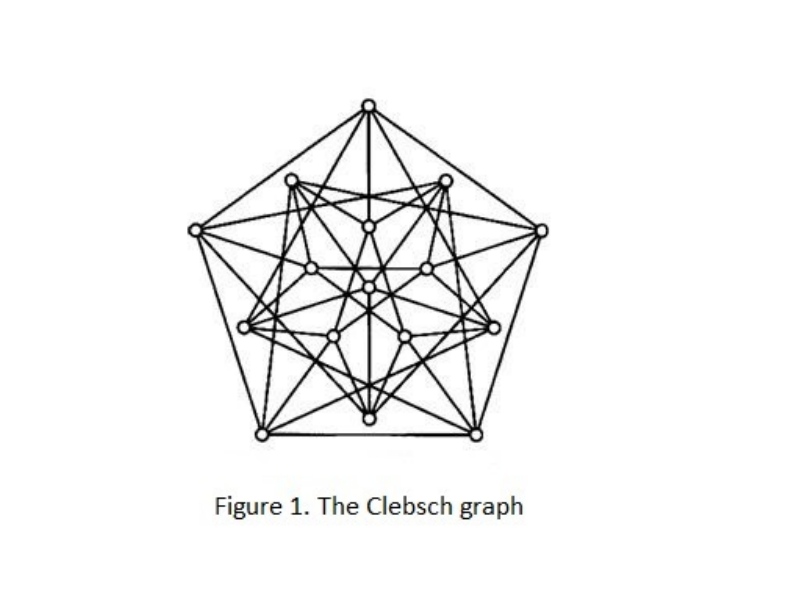}} 
\end{figure}
\\
In this paper,  we are  interested in determining the distance spectrum of the bipartite double cover of  strongly regular graphs. In other words,  let $G=(V,E)$ be a strongly regular graph with parameters $(n,k,a,c)$. We show that there is a close relationship between the spectrum of $G$ and the distance spectrum of $B(G)$. We explicitly determine the distance spectrum of the graph $B(G)$,    according to the spectrum of $G$. Since the spectrum of a strongly regular graph $G$ is determined according to its parameters [8], then the determined   distance spectrum of the graph $B(G)$ by this paper is according to the parameters of the graph $G$.
\section{Main results}
Let $G=(V,E)$ be a graph with an adjacency matrix $A$. In the first step, we show that there is a form for the adjacency matrix   of  the graph $B(G)$, the bipartite double  cover of $G$, according to the matrix $A.$ In the sequel,   $O=O_{n \times n}$ is the matrix in which all entries are zero.
\begin{prop}
Let $G=(V,E)$ be a graph with the vertex-set $V=\{ v_1,v_2,\dots,v_n \}$ and the adjacency matrix $A=(a_{ij})$ in which the rows and columns are indexed by the set $V$. Then $B(G)$, the double cover of $G$, has an adjacency matrix of the form $$ M=\left( \begin{array} {cc}
 O & A \\  A  & O \end{array} \right).  $$
 
\end{prop}
\begin{proof}
Let $V=\{v_1,v_2, \dots, v_n  \}$ be the vertex-set of $G$. Let $P_i=\{ (v_1,i),\dots,$
$(v_n,i) \}$, $i \in \{ 0,1 \}$. Thus $W=P_0 \cup P_1$ is the vertex-set of $B=B(G)$. We know that the rows and columns of $A$ are indexed by the set $V$ such that $A(v_i,v_j)=a_{ij}=1$ if and only if $v_i$ is adjacent to $v_j$. We now index the rows and columns of an adjacency matrix $M =(m_{ij})$ of the graph $B(G)$ by the vertex set $W$ in such a way that   $(v_i,0)$ is adjacent 
to $(v_j,1)$ if and only if $v_i$ is adjacent $v_j$. Hence,  $m_{ij}=1$ if and only if $a_{ij}=1$. Now the result follows.
\end{proof}
Let $G=(V,E)$, $V=\{v_1,v_2, \dots, v_n   \}$ be a connected graph with diameter $d$.
For
every integer $i$,  $0\leq i \leq d$, the distance-$i$ matrix $A_i$ of $G$ is defined as, 
\begin{center}
$A_{i}(v_r,v_s)=\begin{cases}1 & \ if \ d(v_r,v_s)=i \\ 0 & \ otherwise.  \end{cases} $
\end{center}
Then $A_0 = I$ and $A_1$ is the usual adjacency matrix $A$ of $G$. Note   that $A_0 + A_1 + \dots + A_d = J$, where $J$ is the $n \times n$ matrix in which each entry is  1.  Now it is clear that if $D=D(G)$ is the distance matrix of $G$, then 
$$D = A_1 + 2A_2  + 3A_3 + \dots + dA_d.$$
 In the sequel, $J_m$=$J_{m \times m}$ is the all 1 matrix  and $I_m$=$I_{m \times m}$ is the identity matrix of size $m$. \

  In the rest of the paper, we want to focus on $k$-regular graphs of diameter 2. It is easy to check that if $k=2$, then there are only two 2-regular graphs $C_4$  and $C_5$   of diameter 2. Hence, in the sequel we assume that $k \geq 3.$\\
   A graph $G=(V,E)$ is said to be $irreducible$ if for every pair of distinct  vertices $v$ and $w$  in $G$  we have $N(x) \neq N(w)$,  where $N(v)$ denotes the
set of neighbors of the vertex $v$ in  $G$ [17,25]. In other words, the graph $G$ is irreducible if for every pair of vertices $v$ and $w$ when $v \neq w$, then there is a vertex $u$ in $G$ such that $u$ is adjacent to $v$ ($w$) but $u$
  is not adjacent to $w$ ($v$). For instance, the cycle $C_n$, $n \neq 4$ is an irreducible graph but the complete bipartite graph $K_{n,n}, n\geq 2$ is not irreducible. If the graph $G$ is not irreducible, we say that it is $reducible.$
  % Lemma 2.2%
  \begin{lem} 
 Let $G=(V,E)$ be a $k$-regular irreducible graph of order $n$ with  diameter $2$ in which each pair of adjacent vertices have no common neighbor. Let $A$  be an adjacency matrix of 
 $G$.   Let $D$ be the distance matrix of $B(G)$, where $B(G)$ is the bipartite double cover of $G$. Then $D$ can be written  in the following form
 $$D=-2M+2M_4+2X+3Y+2M_5-2I_{2n},$$
  where 
\begin{center}
$   M=\left( \begin{array} {cc}
 O & A \\ A  & O \end{array} \right), \ M_4=\left( \begin{array} {cc}
 A & O \\ O  & A \end{array} \right), \  X=\left( \begin{array} {cc}
 J_n & O \\ O  & J_n \end{array} \right), \  Y=\left( \begin{array} {cc}
 O & J_n \\ J_n  & O \end{array} \right) $  and 
\end{center}  
 $$ M_5=\left( \begin{array} {cc}
 O & I_n \\ I_n  & O \end{array} \right).   $$ 
\end{lem} 
 \begin{proof}
  Let $V=\{v_1,v_2,\dots, v_n  \}$ and $A$ be the adjacency matrix of the graph $G$ in which the  rows and columns are indexed by the set $V.$  Let $P_i=\{ (v_1,i),\dots,$
$(v_n,i) \}$, $i \in \{ 0,1 \}$. Thus $W=P_0 \cup P_1$ is the vertex-set of $B=B(G)$,  the bipartite double cover of  $G$. Hence, by Proposition 2.1, $B$ has an 
adjacency matrix of the form, 
$$ M=\left( \begin{array} {cc}
 O & A \\  A  & O \end{array} \right).  $$ 
Consider the distance-$i$ matrices $M_i$ of the graph $B(G)$. We will show that the diameter of $B(G)$ is 5.  Hence we have 
$$D= M_1 + 2M_2  + 3M_3+4M_4+5M_5,$$
 where  $M_1=M=M_{2n \times 2n}.$ 
We know that in a connected bipartite graph  for every pair of vertices,   the distance between them  is an even integer if and only if they  are in the same part  of its bipartition. Let $x=(v_i,r)$ and $y=(v_j,s)$, $r,s \in \{0,1\} $ be a 
pair of distinct  vertices in $B (G)$.  
 When the diameter of $B(G)$ is 5,  then $ u,v$ are in the same part of $B(G)$ if and only if $d(x,y) \in \{ 2,4 \}$. Since the diameter of $G$ is 2, then each pair of non adjacent vertices of $G$ have at 
least one common neighbor.  Hence $d(x,y)$=2 if and only if $r=s$ and
 $v_i$ and $v_j$ are not adjacent in $G$. Let $r \neq s$ and 
$v_i$  and $v_j$ are not adjacent in $G$ and $v_i \neq v_j$. Since $G$ is irreducible then we can assume that  there is a vertex $u$ in $G$ such that $u$ is adjacent to $v_i$ and $u$ is not adjacent    to $v_j$.  Hence in the graph $B(G)$ the vertex $(v_i,r)$ is adjacent $(u,s)$ and the vertices $(u,s)$ and  $(v_j,s)$ are at distance 2 from each other. Thus,  in $B(G)$  the vertices $(v_i,r)$ and $(v_j,s)$ are at distance 3 from each other. Now assume that $v_i$ and $v_j$ are adjacent.
Since $G$ is a triangle free graph, then $v_i$ and $v_j$ have not a common
neighbor and hence $d((v_i,r),(v_j,r))$$\neq 2$. On the other hand, since $G$ is an irreducible graph, then we can assume that there is a vertex $x$ in $G$ such that $x$ is adjacent to  $v_i$ but $x$ is not adjacent to $v_j$.
Hence, $d((v_j,s),(x,s))=3$ and thus $d((v_i,r),(v_j,s))=4.$\\
We now show that $d((v_i,r),(v_i,s))=5$. Let $x$ be a vertex adjacent to $v_i$. Then, by what is showed up to now,  $d((v_i,s),(x,s))=4$, and hence we have $d((v_i,r),(v_i,s))=5$.
We summarize the
argument in the following array.
\begin{center}
$d(x,y)=\begin{cases} 1 & \   if    \    r \neq s, \   $and$ \    v_i,  v_j    $ are adjacent  in$  \  G,   \ (v_i \neq v_s)    \\ 2 & \ if \ r=s, \  v_i,  v_j   $  are not adjacent  in$  \  G     , \\ 3 &  \   if  \    r \neq s \  $and$ \ v_i,v_j  $  are not   adjacent  in   G$, \\
4 & \  if  \  r=s,  \    v_i,  v_j   $  are  adjacent  in   $  \  G. \\
5  & \  if \  \ r \neq s, \ v_ i=v_j.   
  \end{cases} $
\end{center}
 Now, it is easy to check that 
 \begin{center}
$M_2=\left( \begin{array} {cc}
 J_n-I_n-A & O \\  O  & J_n-I_n-A \end{array} \right)=   \left( \begin{array} {cc}
 J_n & O \\  O  & J_n  \end{array} \right)-\left( \begin{array} {cc}
I_n & O \\  O  & I_n \end{array} \right)-\left( \begin{array} {cc}
A & O \\  O  & A \end{array} \right), \ \ \ \  \ \ \ \ \ \ \ \ \   (1)  $
 
\end{center} 
 
 and 
 
 \begin{center}

$M_3=\left( \begin{array} {cc}
 O & J_n -A-I_n  \\ J_n  -A-I_n  & O \end{array} \right)=\left( \begin{array} {cc}
 O & J_n   \\ J_n   &  O \end{array} \right)-\left( \begin{array} {cc}
O & A \\  A  & O \end{array} \right)-\left( \begin{array} {cc}
O & I_n \\  I_n  & O \end{array} \right),  \ \ \ \  \ \ \ \ \ \  \ \      (2)$  
 
\end{center} 
  \begin{center}

$\ \ \ \ \ \ \ \  M_4=\left( \begin{array} {cc}
 A & O \\ O  & A \end{array} \right),\  M_5=\left( \begin{array} {cc}
O & I_n \\  I_n  & O \end{array} \right). \ \ \  \ \ \ \ \  (3)$  
\end{center} 
If we let $X=\left( \begin{array} {cc}
J_n & O \\  O  & J_n \end{array} \right), Y=\left( \begin{array} {cc}
O & J_n \\  J_n  & O \end{array} \right), $  then we have, \\ 
  $M_2=X-I_{2n}-M_4$, $M_3=Y-M-M_5.$ We now have,
  $$D=M+2M_2+3M_3+4M_4+5M_5=M+2(X-I_{2n} -M_4)+3(Y-M-M_5)+4M_4+5M_5.$$
   Hence, we have,
 $$D=-2M+2M_4+2X+3Y+2M_5-2I_{2n}. \ \ \ \ \ \ \  \ \  (4)$$
\end{proof}
 % Theorem 2.3
 \begin{thm} 
 Let $k \geq 3$ and $G=(V,E)$ be a $k$-regular irreducible graph of order $n$
 with  diameter $2$ in which each pair of adjacent vertices has no common neighbor. Let $A$ be the adjacency matrix of $G$ and  $D$ be the distance matrix of $B(G)$, where $B(G)$ is the bipartite double cover of $G$.  Let the spectrum of $G$ be

  $$Spec(G)=\{k^1,({\lambda_1})^{m_1}, ({\lambda_2})^{m_2}, \dots, ({\lambda_t})^{m_t} \}. $$
  then  the spectrum of $D$ is as follows,\\ 
  $$ Spec(D)= \{ {(5n)}^1, {(4\lambda_1 -4)}^{m_1}, \dots,  {(4\lambda_t -4)}^{m_t}, 0^{n-1}, {(4k-n-4)}^1 \}. $$
  \end{thm}
  \begin{proof} By Lemma 2.2, we know that $D=-2M+2M_4+2X+3Y+2M_5-2I_{2n}$, where $M,M_4,X,Y,M_5$ are the defined matrices in the proof of this lemma.
   Let $j$ be a column of the matrix $J_n$. It is clear that $J_nj=nj.$ Since the rank of $J_n$ is 1, hence we have $Spec(J_n)=  \{ { n}^{1}, 0^{n-1} \}$. Since $G$ is a $k$ regular graph, hence $AJ_n=J_nA=kJ_n.$ Now since $A$ and $J_n$ are symmetric matrices on the field of real numbers $\mathbb{R}$, then there is a basis $B_1=\{ w_1,w_2,\dots,w_n  \}$ for ${\mathbb{R}}^n$
such that each $w_i$ is an eigenvector for both $A$ and $J_n$. We can assume   that $w_1$ is one  for which we have  $J_nw_1=nw_1.$ In fact, we can assume that $w_1=j$.
Note that if $w_1={(x_1,x_2,\dots,x_n)}^t$, then $w_1$ must be in the   eigenspace corresponding to the eigenvalue $n$ of $J_n$. The dimension of this subspace is 1. Hence $w_1=aj$, for some $a \in \mathbb{R}$. Hence we can assume that  $B=\{ w_1=j,w_2,\dots,w_n  \}$. Thus $Aw_1=Aj=kj$. If we let
 $$e_i=\left( \begin{array} {cc}
w_i  \\  w_i  \end{array} \right) \  and \   f_i=\left( \begin{array} {cc}
w_i  \\ - w_i  \end{array} \right), \  1 \leq i \leq n, \ \ \ \ \ \ \ \ \ (5) $$ 
   then we have, \\\\
$Me_i=\left( \begin{array} {cc}
O & A \\  A  & O \end{array} \right)\left( \begin{array} {cc}
w_i  \\  w_i  \end{array} \right)=\left( \begin{array} {cc}
Aw_i  \\  Aw_i  \end{array} \right)=\left( \begin{array} {cc}
\lambda_iw_i  \\  \lambda_iw_i  \end{array} \right)=\lambda_i \left( \begin{array} {cc}
w_i  \\  w_i  \end{array} \right)= \lambda_ie_i$.\\
$Mf_i=\left( \begin{array} {cc}
O & A \\  A  & O \end{array} \right)\left( \begin{array} {cc}
w_i  \\ -w_i  \end{array} \right)=\left( \begin{array} {cc}
-Aw_i  \\  Aw_i  \end{array} \right)=\left( \begin{array} {cc}
-\lambda_iw_i  \\  \lambda_iw_i  \end{array} \right)=-\lambda_i \left( \begin{array} {cc}
w_i  \\ - w_i  \end{array} \right)= -\lambda_if_i$.\\
$M_4e_i=\left( \begin{array} {cc}
A & O \\  O  & A \end{array} \right)\left( \begin{array} {cc}
w_i  \\  w_i  \end{array} \right)=\left( \begin{array} {cc}
Aw_i  \\  Aw_i  \end{array} \right)=\left( \begin{array} {cc}
\lambda_iw_i  \\  \lambda_iw_i  \end{array} \right)=\lambda_i \left( \begin{array} {cc}
w_i  \\  w_i  \end{array} \right)= \lambda_ie_i$.\\
$M_4f_i=\left( \begin{array} {cc}
A & O \\  O  & A \end{array} \right)\left( \begin{array} {cc}
w_i  \\  -w_i  \end{array} \right)=\left( \begin{array} {cc}
Aw_i  \\ - Aw_i  \end{array} \right)=\left( \begin{array} {cc}
\lambda_iw_i  \\  -\lambda_iw_i  \end{array} \right)=\lambda_i \left( \begin{array} {cc}
w_i  \\ - w_i  \end{array} \right)= \lambda_if_i$.\\
$Xe_i=\left( \begin{array} {cc}
J_n & O \\  O  & J_n \end{array} \right)\left( \begin{array} {cc}
w_i  \\  w_i  \end{array} \right)=\left( \begin{array} {cc}
J_n w_i  \\ J_n w_i  \end{array} \right)$. \\ Hence $Xe_1=ne_1$ and if $1<i \leq n$, then $Xe_i=0=0e_i.$ \\ \\
$X f_i=\left( \begin{array} {cc}
J_n & O \\  O  & J_n \end{array} \right)\left( \begin{array} {cc}
w_i  \\  -w_i  \end{array} \right)=\left( \begin{array} {cc}
J_n w_i  \\ - J_n w_i  \end{array} \right)$. \\  Hence $Xf_1=n f_1$ and if $1<i \leq n$, then $X f_i=0=0 f_i.$ \\\\ 
$Ye_i=\left( \begin{array} {cc}
O & J_n \\  J_n  & O \end{array} \right)\left( \begin{array} {cc}
w_i  \\  w_i  \end{array} \right)=\left( \begin{array} {cc}
J_n w_i  \\ J_n w_i  \end{array} \right)$. \\ Hence $Ye_1=ne_1$ and if $1<i \leq n$, then $Y e_i=0=0e_i.$ \\ \\
$Y f_i=\left( \begin{array} {cc}
O & J_n \\  J_n  & O \end{array} \right)\left( \begin{array} {cc}
w_i  \\  -w_i  \end{array} \right)=\left( \begin{array} {cc}
-J_n w_i  \\ J_n w_i  \end{array} \right)$. \\ Hence $Y f_1=-n f_1$ and if $1<i \leq n$, then $Y f_i=0=0 f_i.$ \\ \\
$M_5e_i=\left( \begin{array} {cc}
O & I_n \\  I_n  & O \end{array} \right)\left( \begin{array} {cc}
w_i  \\  w_i  \end{array} \right)=\left( \begin{array} {cc}
I_n w_i  \\ I_n w_i  \end{array} \right)=\left( \begin{array} {cc}
 w_i  \\  w_i  \end{array} \right)$. \\ Hence $M_5e_i=e_i$.\\
$M_5f_i=\left( \begin{array} {cc}
O & I_n \\  I_n  & O \end{array} \right)\left( \begin{array} {cc}
w_i  \\ - w_i  \end{array} \right)=\left( \begin{array} {cc}
-I_n w_i  \\ I_n w_i  \end{array} \right)=\left( \begin{array} {cc}
- w_i  \\  w_i  \end{array} \right)$. \\ Hence $M_5f_i=-f_i$. \\\\
Since $G$ is an irreducible graph, then by Lemma 2.2, we have $D=-2M+2M_4+2X+3Y+2M_5 -2I_{2n}$.   We now have,\\

(i)  $D(e_1)=(-2M+2M_4+2X+3Y+2M_5 -2I_{2n})e_1$=$(-2k+2k+2n+3n+2-2)e_1$=$5ne_1$
=$\mu_1e_1$, where $\mu_1=5n$.\\ \\
 $D(f_1)=(-2M+2M_4+2X+3Y+2M_5 -2I_{2n})f_1$=$(2k+2k+2n-3n-2-2)f_1$=$(4k-n-4)f_1$=$\delta_1 f_1$, where $\delta_1=4k-n-4.$ \\\\
Also for $1 < i \leq n$ we have,\\\\
(ii) 
 $D(e_i)=(-2M+2M_4+2X+3Y+2M_5 -2I_{2n})e_i$=$(-2\lambda_i +2\lambda_i+0+0+2-2)e_i$
 =$0e_i$=$\mu_i e_i$, where $\mu_i=0$. \\
\\
 $D(f_i)=(-2M+2M_4+2X+3Y+2M_5 -2I_{2n})f_i$=$(2\lambda_i+2\lambda_i+0+0-2-2)f_i$
 =$(4\lambda_i -4)f_i$=$\delta_i f_i$, where $\delta_i=4\lambda_i -4$. \\ \\
Noting that $B_1=\{ w_1,w_2, \dots,w_n  \}$ is a basis for $\mathbb{R}^n $, it is easy to check that $B_2=$ $\{ e_1,e_2, \dots, e_n$, $f_1,f_2,\dots, f_n \}$ is a basis of ${\mathbb{R}}^{2n}. $ 
We now conclude the result, that is, 
$$ Spec(D)= \{ {(5n)}^1, {(4\lambda_1 -4)}^{m_1}, \dots,  {(4\lambda_t -4)}^{m_t}, {0}^{n-1},   {(4k-n-4)}^1 \}. $$
\end{proof}
We now consider $k$-regular  graphs of diameter $2$ in which each pair of adjacent vertices has at least one common neighbor.
%  2.4
 \begin{thm} 
 Let $k \geq 3$ and $G=(V,E)$ be a $k$-regular  graph of order $n$ with  diameter $2$ in which each pair of adjacent vertices have at least one common neighbor. Let $A$ be the adjacency matrix of $G$ and $D$ be the distance matrix of $B(G)$, where $B(G)$ is the bipartite double cover of $G$.  Let the spectrum of $G$ be

  $$Spec(G)=\{k^1,({\lambda_1})^{m_1}, ({\lambda_2})^{m_2}, \dots, ({\lambda_t})^{m_t} \}. $$
  then  the spectrum of $D$ is \\ 
  $$ Spec(D)= 
 \{ {(-2k+5n-2)}^1, {(2\lambda_1 -2)}^{m_1}, \dots,  {(2\lambda_t -2)}^{m_t},$$
$$ {(-2\lambda_t-2)}^{m_t}, \dots, {(-2\lambda_1-2)}^{m_1}, {(2k-n-2)}^1 \}. $$
  \end{thm}
  \begin{proof}
Let $V=\{v_1,v_2,\dots, v_n  \}$ and $A$ be the adjacency matrix of the graph $G$ in which the  rows and columns are indexed by the set $V$. Let $P_i=\{ (v_1,i),\dots,$
$(v_n,i) \}$, $i \in \{ 0,1 \}$. Thus $W=P_0 \cup P_1$ is the vertex set of $B=B(G)$,  the bipartite double cover of  $G$. Thus by Proposition 2.1, $B(G)$ has an 
adjacency matrix of the form, 
$$ M=\left( \begin{array} {cc}
 O & A \\  A  & O \end{array} \right).  $$ 
 Let $x=(v_i,r), y=(v_j,s)$ be two distinct vertices in the graph $B=B(G)$. Since $v_i$ and $v_j$ have at least one common neighbor in the graph $G$, thus if $r=s$, that is, $x$ and $y$ are in the same part of the bipartite graph $B$, then they are at distance 2 in $B$. Moreover if $x$ and $y$ are not adjacent and are not in the same part of $B$, then they are at distance 3 in this graph. To check this, let $z=(u,s)$ be an adjacent vertex to $x=(v_i,r)$ in $B(G)$. Now, since $d(z,y)=2$, then we have $d(x,y)=3.$      Hence the diameter of $B(G)$ is 3.
Now, by an argument similar to what we have done  in the proof  of Lemma 2.2, we deduce that
$$ D=\left( \begin{array} {cc}
O & A \\  A  & O \end{array} \right)+2\left( \begin{array} {cc}
J_n-I_n & O \\  O  & J_n-I_n \end{array} \right)+3\left( \begin{array} {cc}
O& J_n-A \\  J_n-A  & O \end{array} \right).  $$
Thus,
$$ D=-2M+2X+3Y -2I_{2n},  \ \ \ \ \ \ \ \ \ (6)$$ \\
where $X$ and $Y$ are the matrices which have been defined in the proof of Lemma 2.2. Since $A$ is a $k$ regular graph, then $AJ_n=J_nA=kJ$. It is easy to see that the set $S=\{ M,X,Y \}$ is a commuting set 
of real symmetric matrices. Let $j$ be a column of the matrix $J_n$. Now by an argument similar to what we did in the proof of Theorem 2.3, we deduce that 
there is a basis $B_1=\{ w_1=j,w_2, \dots,w_n  \}$ for ${\mathbb{R}}^n$ such that each element of $B_1$ is an  eigenvector for $A$ and $J_n$. Let $e_i$ and $f_i$ be the column matrices which are defined in (5).
Now, by a similar argument which we have done in the proof of Theorem 2.3, we deduce that \\\\
 (i)  $D(e_1)=(-2M+2X+3Y -2I_{2n})e_1$=$(-2k +2n+3n-2)e_1$=$(-2k+5n-2)e_1$
=$\mu_1e_1$, where $\mu_1=-2k+5n-2$.\\ \\
 $D(f_1)=(-2M+ 2X+3Y -2I_{2n})f_1$=$(2k+2n-3n-2)f_1$=$(2k-n-2)f_1$=$\delta_1 f_1$, where $\delta_1=2k-n-2.$ \\\\
Also for $1 < i \leq n$ we have,\\\\
(ii) 
 $D(e_i)=(-2M+2X+3Y -2I_{2n})e_i$=$(-2\lambda_i+0+0-2)e_i$
 =$(-2\lambda_i-2)e_i$=$\mu_i e_i$, where $\mu_i=-2\lambda_i-2$. \\
\\
 $D(f_i)=(-2M +2X+3Y -2I_{2n})f_i$=$(2\lambda_i +0+0-2)f_i$
  =$\delta_i f_i$, where $\delta_i=2\lambda_i -2$. \\ \\
Noting that $B_1=\{ w_1,w_2, \dots,w_n  \}$ is a basis for $\mathbb{R}^n $, it is easy to check that $B_2=$ $\{ e_1,e_2, \dots, e_n$, $f_1,f_2,\dots, f_n \}$ is a basis of ${\mathbb{R}}^{2n}. $ We now conclude the result, that is, 
$$ Spec(D)= 
 \{ {(-2k+5n-2)}^1, {(2\lambda_1 -2)}^{m_1}, \dots,  {(2\lambda_t -2)}^{m_t},$$
$$ {(-2\lambda_t-2)}^{m_t}, \dots, {(-2\lambda_1-2)}^{m_1}, {(2k-n-2)}^1 \}. $$
\end{proof}
 Let $G$ be a strongly regular graph of parameters $(n,d,a,c)$. Then $G$ has three eigenvalues $d$, $\lambda_1$ and $ \lambda_2$ where $\lambda_1 =\frac {(a-c)+\sqrt {\Delta}}{2}$, $\Delta=(a-c)^2+4(d-c)$ and $\lambda_2=\frac {(a-c)-\sqrt {\Delta}}{2} $  [8]. It is clear that the multiplicity of $d$ is 1. If $m_{\lambda_1}$ and $m_{\lambda_2}$ are multiplicities of $\lambda_1$ and $\lambda_2$  respectively, then
 $$ m_{\lambda_1}=\frac{1}{2}((n-1)-\frac{2d+(n-1)(a-c)}{\sqrt \Delta})$$
 and                                                                                  
$$ m_{\lambda_2}=\frac{1}{2}((n-1)+\frac{2d+(n-1)(a-c)}{\sqrt \Delta}). \   \ \ \  \ \ \ \ (7)$$ 
 
It is clear that the diameter of a strongly regular graph $G$  is 2. But, it is not true that the diameter of $B(G)$, the bipartite double cover of $G$,  is always  3. By Theorem 2.3,  if $a \neq 0$ then the diameter of $B(G)$ is 3, and if $a = 0$ and $G$ is an irreducible graph, then the diameter of $B(G)$ is 5. There are strongly regular graphs of parameters $(n,d,a,c)$ in which we have $a=0$. For instance, the Petersen graph is strongly regular graph with the parameters $(10,3,0,1)$. 
Also, the Hoffman-Singleton graph  is a strongly regular graph with the parameters $(50,7,0,1)$ [4,8]. For more information about strongly regular graphs with the parameters $(n,d,a,c)$ in which $a=0$  (see [4]). \\
 It is quite possible that a strongly regular graph with parameters $(n,k,0,c)$,  that is a triangle free strongly regular graph,  be   reducible. For instance the complete bipartite graph $K_{m,m}, m \geq 2$ is a triangle free  strongly regular graph of parameter $(2n,n,0,n)$ which is reducible. In the following lemma, we show that this is an exceptional case.
\begin{lem}
Let $G=(V,E)$ be a connected  reducible triangle free strongly regular graph. Then  $G$ is isomorphic with the complete bipartite graph $K_{m,m}$ for some positive integer $m \geq 2$. 
\end{lem}
\begin{proof}
Let $G=(V,E)$ be a connected reducible triangle free strongly regular graph with parameters $(n,m,0,c)$. If $m=1$, then since $G$ is connected we have $G=K_2$  which is irreducible. Hence we assume that $m \geq 2.$ Since $G$ is 
reducible, there is a pair of distinct  vertices $v,w$ in $G$ such that $N(v)=N(w)$. Note that when $v$ and $w$ are adjacent, we have $w \in N(v)$ but $w \notin N(w)$, which implies that $N(w) \neq N(v).$ Hence we deduce that   $v$ and $w$ are not adjacent. Thus, there is a pair of non adjacent vertices $v,w$ in $G$ such that they have $|N(v)|=m$ common neighbors. Now, since $G$ is a strongly regular graph, we deduce that every pair of distinct non adjacent vertices in $G$ have $m$ common neighbors, that is $c=m$. Let $P=N(v)$ and $Q=V-N(v)=V-P.$ It is clear that $v,w \in Q. $ If $x \neq v$ is a vertex of $G$ in $Q$,  then $x$ and $v$ are not adjacent, hence they have $m$ common neighbors, which implies that $N(x)=N(v). $ In oder words, each vertex in $Q$ 
is adjacent to every vertex in $P$. Nothing that $G$ is an $m$-regular graph, we deduce that each pair of distinct vertices in $Q$ are non adjacent. On the other hand, since $G$ is triangle free graph, then each pair of distinct vertices in $P$ must be non adjacent. Therefore, $G$ is an $m$-regular bipartite graph. Hence $|P|=|Q|=m$. We now conclude that $G$ is isomorphic with the  complete  bipartite  graph $K_{m,m}.$ 
\end{proof}
We now,  by Theorem 2.3  and Theorem 2.4 and Lemma 2.5,   can determine the distance spectrum of bipartite double cover of strongly regular graphs.   
 \begin{thm}
 Let $G=(V,E)$ be a strongly regular graph with parameters $(n,d,a,c)$ and the spectrum $\{ d^1, {\lambda_1}^{m_1}, {\lambda_2}^{m_2} \}.$ let $D$ be the distance matrix of the graph $B(G)$, the bipartite double cover of $G$. If $a \neq 0$,  then we have,
  $$Spec(D) =$$
  $$ \{ {(-2d+5n-2)}^1, {(2\lambda_1-2)}^{m_1}, {(2\lambda_2-2)}^{m_2},{(-2\lambda_2-2)}^{m_2},{(-2\lambda_1-2)}^{m_1},$$
  $${(2d-n-2)}^1 \} $$
 and if $a=0$ and $G \ncong K_{m,m}$, $ m \geq 2, $ then we have,
 $$Spec(D) = 
\{ {(5n)}^1, {(4\lambda_1-4)}^{m_1}, {(4\lambda_2-4)}^{m_2}, {0}^{n-1}, {(4d-n-4)}^1 \}.  $$
\end{thm} 
 \section{Some examples}
 
 (i) The Petersen graph $P$ is a strongly regular graph with the parameters $(10,3,0,1)$ [8]. By $(7)$, we can check that, $Spec(G)= \{ 3^1, 1^5,{(-2)}^{4} \}$. Let $D$ be the distance matrix of the graph $B(P)$, the bipartite double cover of $G$. Hence by Theorem 2.6, we can check that
  $$Spec(D)=\{ {(50)}^1,0^{14},{(-12)}^{4}, {-2}^{1} \}.$$ 
(ii) As we stated, the Hoffman-Singleton graph  is a  strongly regular graph with the parameters $(50,7,0,1)$ [8]. Hence, by $(7)$, we can check that, $Spec(G)= \{ {7}^1, 2^{28},{(-3)}^{21} \}$. Let $D$ be the distance matrix of the graph $B(G)$, the bipartite double cover of $G$. Thus by Theorem 2.6, we can check that 
$$Spec(D)=\{ {(250)}^1,4^{28},{(-16)}^{21}, {0}^{49},{(-26)}^1 \}.$$ 
 (iii)  The line graph $L(K_{5,5)}=G$ is a strongly regular graph with parameters $(25,8,3,2)$ [8]. Hence, by $(7)$, we can check that, $Spec(G)= \{ 8^1, 3^8,{(-2)}^{16} \}$. Let $D$ be the distance matrix of the graph $B(G)$, the bipartite double cover of $G$. Thus by Theorem 2.6, we can check that
  $$Spec(D)=\{ {(107)}^1,4^8,{(-6)}^{16}, 2^{16},{(-8)}^{16}, {(-11)}^1 \}.$$
(iv) Let $n$ be a positive integer and $\Gamma$ be a group of order $n$ with the identity element 1. Consider the group $\Gamma_1=\Gamma \times \Gamma.$  Let $S= \{ (g,1),(1,g),(g,g) | 1\neq g \in \Gamma  \}$. It is not difficult to check that the Cayley graph $G=Cay(\Gamma_1,S)$ is a strongly regular graph with the parameters $(n^2,3n-3,n,6)$ [23]. By $(7)$, we can check that, \\
   $Spec(G)= \{ {(3n-3)}^1, {(n-3)}^{m_1},{(-3)}^{m_2} \}$, where $m_1=\frac{1}{2}(n^2 -n)$,  $m_2=\frac{1}{2}(n^2+n-2). $
Thus by Theorem 2.6, we can check that 
 Spec(D)=$$\{ {(5n^2-6n+4)}^1,{(2n-8)}^{m_1},{(-8)}^{m_2}, {4}^{m_2},{(-2n+4)}^{m_1},{(-n^2+6n-8)^1} \}.$$  
\section{Conclusion}
In this paper, we have determined the distance spectrum of the bipartite double cover of strongly regular graphs according to their parameters  (Theorem 2.5). Also, we have determined the distance spectrum of bipartite double cover of  some other 
 classes of graphs with diameter 2 according to their spectrum (Theorem 2.3 and Theorem 2.4). In all the discussed cases, we saw that if the strongly regular graph $G=(V,E)$ is integral, then it is distance integral.
 
\section{ Declarations   }
    
{\bf Conflicts of interest }
\\\\
 The corresponding author states that there is no conflict of
interest.


\begin{thebibliography}{}
%
\bibitem{1}  Aouchiche M,  Hansen P.  Distance spectra of graphs: a survey,  Linear Algebra Appl.
458 (2014), 301-386.
%
\bibitem{2}  Bali$\acute{a}$ska  K,     Cvetkovi$\acute{c}$ D,   Radosavljevi$\acute{c}$ Z,    Simi$\acute{c}$ S,  Stevanovi$\acute{c}$ D.   A survey on integral
graphs,  Publ. Elektroteh. Fak., Univ. Beogr., Ser. Mat. 13 (2002), 42–65.
%

\bibitem{3}  Biggs  N.  Strongly regular graphs with no triangles,  arXiv preprint arXiv:0911.2160, (2009).
% 
\bibitem{4}   Brouwer A.E,    Cohen A.M,   Neumaier A. Distance-Regular Graphs, Springer-
Verlag, New York, (1989).
%


\bibitem{5} Chai Y,  Wang L,   Zhou  Y.   DQ-integral and DL-integral generalized wheel graphs,   Indian Journal of Pure and Applied Mathematics (2025): 1-12.
%
\bibitem{6} Cvetkovic  D, Rowlinson  P,  Simic  S.   An introduction to the theory of graph spectra,   Cambridge University Press, (2010).
%
\bibitem{7} Donno A.  "Spectrum, distance spectrum, and Wiener index of wreath products of complete graphs.",  Ars Mathematica Contemporanea 13, no. 1 (2017): 207-225.
%
\bibitem{8} Godsil C,  Royle  G.   Algebraic Graph Theory,  Springer, (2001).
%
\bibitem{9}    Harary F,    Schwenk A.J. Which graphs have integral spectra?, In Graphs and Combinatorics, (eds. R. Bari and F. Harary), (Proc. Capital Conf., George Washington Univ., Washington, D.C., (1973), Lecture Notes in Mathematics
406, Springer-Verlag, Berlin (1974), 45-51.
%
\bibitem{10}  Huang J,   Li  S.C.   Distance integral Cayley graphs over abelian groups and dicyclic
groups, J. Algebraic Combin. 53 (2021) 921-943.
%
\bibitem{11}  Indulal G. "Distance spectrum of graph compositions.",  Ars Math. Contemp. 2, no. 1 (2009): 93-100.


%
\bibitem{12} Kogani  R,   Mirafzal  S.M.  On determining the distance spectrum of a class of distance
integral graphs, J. Algebr. Syst, 10(2) (2023), 299–308.
%
\bibitem{13} Lin H,  Shu J,   Xue J,   Zhang  Y.  A survey on distance spectra of graphs, Adv. Math.
(China) 50(1) (2021) 29-76.
%
\bibitem{14} Lipschutz  S,  Abellanas  L,   Ontalba  C.M.  Linear Algebra (Vol. 2). Madrid: McGraw-Hill, (1992).
%
\bibitem{15} Mirafzal S.M.   A new class of  integral  graphs constructed from the hypercube,     Linear Algebra Appl. 558 (2018) 186-194.
%
\bibitem{16} Mirafzal  S.M.    The automorphism group of the bipartite Kneser graph, Proceedings-Mathematical Sciences, (2019), doi.org/10.1007/s12044-019-0477-9.

%
\bibitem{17} Mirafzal  S.M.  On the automorphism groups of connected bipartite irreducible graphs,  Proc. Math.
Sci. (2020). https://doi.org/10.1007/s12044-020-0589-1.
%
\bibitem{18} Mirafzal  S.M. Some remarks on the square graph of the hypercube,  Ars Mathematica  Contemporanea    Vol. 23  No. 2 (2023)
%
\bibitem{19}  Mirafzal  S.M.   The line graph of the crown graph is distance integral,  Linear and Multilinear Algebra 71, no. 4 (2023): 662-672.
%
\bibitem{20} Mirafzal  S.M.  On the distance eigenvalues of design graphs,  Ricerche mat 73, 2759–2769 (2024),  https://doi.org/10.1007/s11587-023-00794-w
%
\bibitem{21} Mirafzal  S.M. The distance spectrum of the line graph of the
crown graph, arXiv:2508.07202v1.
%


\bibitem{22} Mirafzal  S.M. The distance spectrum of regular bipartite graphs of  diameter 3,  Ricerche di Matematica (2026),  
https://doi.org/10.1007/s11587-026-01078-9.

\bibitem{23} Nica   B.   A Brief Introduction to Spectral Graph Theory, EMS Publishing House, Zuerich, (2018). 
%
\bibitem{24}  Pokorn$\acute{y}$ M,  H$\acute{i}$c P,   Stevanovi$\acute{c}$ D,  Milo$\breve{s}$evi$\acute{c}$ M.  On distance integral graphs,
Discrete Math, 338 (2015), 1784–1792.
%
\bibitem{25} Sabidussi G. Vertex-transitive graphs,   Monatsh. Math. 68. 426–438 (1964).
%
\bibitem{26}Wu Y,  Zhang X,  Feng L, Wu T.  Distance and adjacency spectra and eigenspaces for three (di) graph lifts: A unified approach,  Linear Algebra and its Applications. (2023)  Sep 1;672:147-81.
%
\bibitem{27} Zhang Y, Lin H. Perfect matching and distance spectral radius in graphs and
bipartite graphs. Discrete Applied Mathematics (2021), 304, pp.315-322.
%
\bibitem{28}  Zou  L, Wu  Y,    Feng  L.  DL-integral and DQ-integral $n$-Cayley graphs. Applied Mathematics and Computation (2025), 489, p.129178. 
%
\end{thebibliography}
\end{document}